\documentclass{amsart}

\theoremstyle{definition} \theoremstyle{example}
\theoremstyle{remark} \numberwithin{equation}{section}

\begin{document}

\centerline{INFINITE ORDER DECOMPOSITIONS OF C$^*$-ALGEBRAS}

\bigskip

\centerline{\it Arzikulov Farhodjon Nematjonovich}

\medskip

subjclass: {Primary 54C40, 14E20;      Secondary 46E25, 20C20}

\medskip

{\bf Keywords:} {\it C$^*$-algebra, Peirce decomposition, von
Neumann algebra}

\bigskip

\centerline{\bf Abstract}

\medskip

In the given article infinite order decompositions of
C$^*$-algebras are investigated. We give complete proofs of the
following statements:

1) If the order unit space $\sum_{\xi,\eta}^\oplus p_\xi Ap_\eta$
is monotone complete in $B(H)$ (i.e. ultraweakly closed), then
$\sum_{\xi,\eta}^\oplus p_\xi Ap_\eta$ is a C$^*$-algebra.

2) If $A$ is monotone complete in $B(H)$ (i.e. a von Neumann
algebra), then $A=\sum_{\xi,\eta}^\oplus p_\xi Ap_\eta$.

3) If $\sum_{\xi,\eta}^\oplus p_\xi Ap_\eta$ is a C$^*$-algebra
then this algebra is a von Neumann algebra.

\bigskip

\section*{ Introduction}

\medskip

In the given article the notion of infinite order decomposition of
a C*-algebra with respect to an infinite orthogonal set of
projections is investigated. It is known that for any projection
$p$ of a C$^*$-algebra $A$ the next equality is valid $A=pAp\oplus
pA(1-p)\oplus (1-p)Ap \oplus (1-p)A(1-p)$, where $\oplus$ is a
direct sum of spaces. In the given article we investigate infinite
analog of this decomposition, namely, infinite order
decomposition. The notion of infinite order decomposition was
introduced in \cite{1}. The following theorems belong to \cite{1}:

let $A$ be a C$^*$-algebra on a Hilbert space $H$, $\{p_\xi\}$ be
an infinite orthogonal set of projections in $A$ with the least
upper bound $1$ in the algebra $B(H)$. Then

{\it 1) if the order unit space $\sum_{\xi,\eta}^\oplus p_\xi
Ap_\eta$ is monotone complete in $B(H)$ (i.e. ultraweakly closed),
then $\sum_{\xi,\eta}^\oplus p_\xi Ap_\eta$ is a C$^*$-algebra,

2) if $A$ is monotone complete in $B(H)$ (i.e. a von Neumann
algebra), then $A=\sum_{\xi,\eta}^\oplus p_\xi Ap_\eta$,

3) if $\sum_{\xi,\eta}^\oplus p_\xi Ap_\eta$ is a C$^*$-algebra
then this algebra is a von Neumann algebra.}

In the given article we give complete proofs of these statements
(see, respectively, item 2) of theorem 11, proposition 9, item 2)
of corollary 14). Also we prove that for the infinite order
decomposition $\sum_{\xi,\eta}^\oplus p_\xi Ap_\eta$ of a
C$^*$-algebra $A$ with respect to an infinite orthogonal set
$\{p_i\}$ of projections of $A$, if $p_\xi Ap_\xi$ is a von
Neumann algebra for any $\xi$ then $\sum_{\xi,\eta}^\oplus p_\xi
Ap_\eta$ is a von Neumann algebra. For this propose it was
constructed operations of multiplication and involution
corresponding to infinite order decompositions. It turns out, the
order and the norm defined in the infinite order decomposition of
a C$^*$-algebra on a Hilbert space $H$ coincide with the usual
order and the norm in the algebra $B(H)$. Also, it is proved that,
if a C$^*$-algebra $A$ with an infinite orthogonal set $\{p_\xi\}$
of projections in $A$ such that $\sup_\xi p_\xi=1$ is not a von
Nemann algebra, projections of the set $\{p_\xi\}$ are pairwise
equivalent then $A\neq \sum_{\xi,\eta}^\oplus p_\xi Ap_\eta$.
Moreover if the order unit space $\sum_{\xi,\eta}^\oplus p_\xi
Ap_\eta$ is not weakly closed then $\sum_{\xi,\eta}^\oplus p_\xi
Ap_\eta$ is not a C$^*$-algebra.

The author wants to thank V.I.Chilin and A.A.Rakhimov for many
stimulating conversations on the subject.

\newpage

\bigskip

\section{Infinite order decompositions}

\medskip

Let $A$ be a C$^*$-algebra on a Hilbert space $H$, $\{p_\xi\}$ be
an infinite orthogonal set of projections of $A$ with the least
upper bound (LUB) $1$, calculated in $B(H)$. By
$\sum_{\xi,\eta}^\oplus p_\xi Ap_\eta$ we denote the set
$$
\{\{a_{\xi,\eta}\}:   a_{\xi,\eta}\in p_\xi Ap_\eta\,\, for\, all
\,\,\xi, \eta, \,and\,there\, exists\, such\, number
$$
$$
K\in R \,\,that \Vert\sum_{k,l=1}^na_{kl}\Vert\le K \,\,for\,
all\,\, n\in N \,\,and\,\,\,\{a_{kl}\}_{kl=1}^n\subseteq
\{a_{\xi,\eta}\}\},
$$
and say $\sum_{\xi,\eta}^\oplus p_\xi Ap_\eta$ is {\it an infinite
order decomposition} of $A$.

We define a relation of order $\leq$ in the vector space
$\sum_{\xi,\eta}^\oplus p_\xi Ap_\eta$ as follows: for elements
$\{a_{\xi\eta}\}$, $\{b_{\xi\eta}\}\in \sum_{\xi,\eta}^\oplus
p_\xi Ap_\eta$, if for all $n\in N$, $\{p_k\}_{k=1}^n\subset
\{p_\xi\}$ the inequality $\sum_{k,l=1}^na_{kl}\leq
\sum_{k,l=1}^nb_{kl}$ holds, then we will write
$\{a_{\xi\eta}\}\leq \{b_{\xi\eta}\}$. Also, the map
$\{a_{\xi,\eta}\}\to\Vert \{a_{\xi,\eta}\}\Vert$,
$\{a_{\xi,\eta}\}\in \sum_{\xi,\eta}^\oplus p_\xi Ap_\eta$, where
$\Vert \{a_{\xi,\eta}\}\Vert=\sup \{\Vert \sum_{kl=1}^n
a_{kl}\Vert :n\in N, \{a_{kl}\}_{kl=1}^n \subseteq
\{a_{\xi,\eta}\}\}$, is a norm on the vector space
$\sum_{\xi,\eta}^\oplus p_\xi Ap_\eta$.

{\it Example.} Let $n$ be an infinite cardinal number, $\Xi$ be a
set of indexes of the cardinality $n$. Let $\{e_{ij}\}$ be a set
of matrix units such that $e_{ij}$ is a $n\times n$-dimensional
matrix, i.e. $e_{ij}=(a_{\alpha\beta})_{\alpha\beta\in\Xi}$, the
$(i,j)$-th component of which is $1$, i.e. $a_{ij}=1$, and the
rest components are zeros. Let $\{m_\xi\}_{\xi\in \Xi}$  be a set
of $n\times n$-dimensional matrices. By $\sum_{\xi\in \Xi} m_\xi$
we denote the matrix whose components are sums of the
corresponding components of matrices of the set $\{m_\xi
\}_{\xi\in \Xi}$. Let
$$
M_n({\bf C})=\{\{\lambda_{ij}e_{ij}\}: \,for\,\, all\,\,
indexes\,\, i,\,j \,\lambda_{ij}\in {\bf C},
$$
$$
and\,\, there\,\, exists\,\, such\,\, number\,\, K\in {\bf
R},\,\,that \,\, for \,\, all\,\, n\in N
$$
$$
and\,\, \{e_{kl}\}_{kl=1}^n\subseteq \{e_{ij}\} \Vert\sum_{kl=1}^n
\lambda_{kl}e_{kl}\Vert \le K\},
$$
where $\Vert \,\, \Vert$ is a norm of a matrix. It is clear that
$M_n({\bf C})$ is a vector space. The set $M_n({\bf C})$, defined
above, coincides with the set
$$
\mathcal{M}_n({\bf C})=\{\{\lambda_{ij}e_{ij}\}:\,for\,\, all\,\,
indexes\, ij\,\,\lambda_{ij}\in {\bf C},
$$
$$
and\,\,\, there\,\,\, exists\,\,\, such\,\,\, number\,\,\, K\in R
\,\,that \,\, for\,\, all\,
$$
$$
\{x_i\}\in l_2(\Xi)\,\, the\,\,\, next\,\,\, inequality\,\,\,
holds \, \sum_{j\in \Xi} \vert\sum_{i\in \Xi}\lambda_{ij}x_i
\vert^2\le K^2\sum_{i\in \Xi} \vert x_i\vert^2\},
$$
where $l_2(\Xi)$ is the Hilbert space on ${\bf C}$ with elements
$\{x_i\}_{i\in \Xi}$, where $x_i\in {\bf C}$ for all $i\in \Xi$.

Associative multiplication of elements in $M_n({\bf C})$ can be
defined as follows: if $x=\sum_{ij\in \Xi}\lambda_{ij}e_{ij}$, $
y=\sum_{ij\in \Xi}\mu_{ij}e_{ij}$ are elements of $M_n({\bf C})$
then
$$
xy=\sum_{ij\in \Xi} \sum_{\xi\in \Xi} \lambda_{i\xi}\mu_{\xi
j}e_{ij}.
$$
With this operation $M_n({\bf C})$ is an associative algebra and
$M_n({\bf C})\equiv B(l_2(\Xi))$, where $B(l_2(\Xi))$ is the
associative algebra of all bounded linear operators on the Hilbert
space $l_2(\Xi)$. Then $M_n({\bf C})$ is a von Neumann algebra of
infinite $n\times n$-dimensional matrices over ${\bf C}$.

Similarly, if we take the algebra $B(H)$ of all bounded linear
operators on a Hilbert space $H$ and if $\{q_i\}$ is a maximal
orthogonal set of minimal projections of the algebra $B(H)$, then
$B(H)=\sum_{ij}^\oplus q_i B(H)q_j$ (see \cite{1}).

\bigskip

Let $A$ be a C$^*$-algebra on a Hilbert space $H$, $\{p_i\}$ be an
infinite orthogonal set of projections with the LUB $1$ in $B(H)$
and $\mathcal{A}=\{ \{p_iap_j\}: a\in A\}$. Then $A\equiv
\mathcal{A}$ \cite{2}.

{\bf Lemma 1.} {\it Let $A$ be a C$^*$-algebra on a Hilbert space
$H$, $\{p_\xi\}$ be an infinite orthogonal set of projections of
$A$ with the LUB $1$ in $B(H)$. Then, $\sum_{\xi,\eta}^\oplus
p_\xi Ap_\eta$ is a vector space with the following componentwise
algebraic operations
$$
\lambda\cdot \{a_{\xi\eta}\}=\{ \lambda a_{\xi\eta}\}, \lambda\in
{\bf C}
$$
$$
\{a_{\xi\eta}\}+\{b_{\xi\eta}\}=\{a_{\xi\eta}+b_{\xi\eta}\},
a_{\xi\eta}, b_{\xi\eta}\in \sum_{\xi,\eta}^\oplus p_\xi Ap_\eta .
$$
And the vector space $\mathcal{A}$ is a vector subspace of the
vector space $\sum_{\xi,\eta}^\oplus p_\xi Ap_\eta$. }

\medskip

{\bf Lemma 2.} {\it Let $A$ be a C$^*$-algebra on a Hilbert space
$H$, $\{p_\xi\}$ be an infinite orthogonal set of projections of
$A$ with the LUB $1$ in $B(H)$. Then, the map
$\{a_{\xi,\eta}\}\to\Vert \{a_{\xi,\eta}\}\Vert$,
$\{a_{\xi,\eta}\}\in \sum_{\xi,\eta}^\oplus p_\xi Ap_\eta$, where
$\Vert \{a_{\xi,\eta}\}\Vert=\sup \{\Vert \sum_{kl=1}^n
a_{kl}\Vert :n\in N, \{a_{kl}\}_{kl=1}^n \subseteq
\{a_{\xi,\eta}\}\}$, is a norm, and $\sum_{\xi,\eta}^\oplus p_\xi
Ap_\eta$ is a Banach space with this norm.}

{\bf Proof.} It is clear, that for any element
$\{a_{\xi,\eta}\}\in \sum_{\xi,\eta}^\oplus p_\xi Ap_\eta$, if
$\Vert \{a_{\xi,\eta}\}\Vert=0$, then $a_{\xi,\eta}=0$ for all
$\xi$, $\eta$, i.e. $\{a_{\xi,\eta}\}=0$. The rest conditions in
the definition of the norm can be also easily checked.

Let $(a_n)$ be a Cauchy sequence in $\sum_{\xi,\eta}^\oplus p_\xi
Ap_\eta$, i.e. for any positive number $\varepsilon >0$ there
exists $n\in {\bf N}$ such, that $\Vert
a_{n_1}-a_{n_2}\Vert<\varepsilon$ for all $n_1\geq n$, $n_2\geq
n$. Then the set $\{\Vert a_n\Vert\}$ is bounded by some number
$K\in {\bf R}_+$ and for any finite set $\{p_k\}_{k=1}^n\subset
\{p_i\}$ the sequence $(pa_np)$ is a Cauchy sequence, where
$p=\sum_{k=1}^n p_k$. Then, since $A$ is a Banach space, we have
$\lim_{n\to \infty} pa_np\in A$.

Let $a_{\xi,\eta}=\lim_{n\to \infty} p_\xi a_np_\eta$ for all
$\xi$ and $\eta$. Then $\Vert \sum_{kl=1}^n a_{kl}\Vert\leq K$ for
all $n\in {\bf N}$ and $\{a_{kl}\}_{kl=1}^n \subseteq
\{a_{\xi,\eta}\}$. Hence $\{a_{\xi,\eta}\}\in
\sum_{\xi,\eta}^\oplus p_\xi Ap_\eta$. $\triangleright$

\medskip

The definition of the order in $\sum_{\xi,\eta}^\oplus p_\xi
Ap_\eta$ is equivalent to the following condition: for the
elements $\{a_{\xi\eta}\}$,
$\{b_{\xi\eta}\}\in\sum_{\xi,\eta}^\oplus p_\xi Ap_\eta$, if
$\{a_{kl}\}_{k,l=1}^n\leq \{b_{kl}\}_{k,l=1}^n$ for all $n\in N$
and $\{p_k\}_{k=1}^n \subseteq \{p_i\}$ in $\mathcal{A}$, then
$\{a_{\xi\eta}\}\leq \{b_{\xi\eta}\}$.

{\bf Proposition 3.} {\it Let $A$ be a C$^*$-algebra on a Hilbert
space $H$, $\{p_\xi\}$ be an infinite orthogonal set of
projections in $A$ with the LUB $1$ in $B(H)$. Then the relation
$\leq$, introduced above, is a relation of a partial order, and
$\sum_{\xi,\eta}^\oplus p_\xi Ap_\eta$ is an order unit space with
this order. In this case $\mathcal{A}=\{\{p_\xi ap_\eta\}: a\in
A\}$ is an order unit subspace of the order unit space
$\sum_{\xi,\eta}^\oplus p_\xi Ap_\eta$. }

{\bf Proof.} Let $\mathcal{M}=\sum_{\xi,\eta}^\oplus p_\xi
Ap_\eta$. $\mathcal{M}$ is a partially ordered vector space, i.e.
$\mathcal{M}_+\cap \mathcal{M}_-=\{0\}$, where $\mathcal{M}_+=\{
\{a_{\xi\eta}\}\in \mathcal{M}: \{a_{\xi\eta}\}\geq 0\}$,
$\mathcal{M}_-=\{ \{a_{\xi\eta}\}\in \mathcal{M}:
\{a_{\xi\eta}\}\leq 0\}$.

By the definition of the order $\mathcal{M}$ is Archimedean. Let
$\{a_{\xi\eta}\}\in \mathcal{M}$. Since $-\Vert
\{a_{\xi,\eta}\}\Vert p\leq p\{a_{\xi,\eta}\}p\leq \Vert
\{a_{\xi,\eta}\}\Vert p$ for any finite set
$\{p_k\}_{k=1}^n\subset \{p_\xi\}$, where $p=\sum_{k=1}^n p_k$, we
have $-\Vert \{a_{\xi,\eta}\}\Vert 1\leq \{a_{\xi,\eta}\}\leq
\Vert \{a_{\xi,\eta}\}\Vert 1$ by the definition of the order, and
the unit of $A$ is an order unit of the partially ordered vector
space $\mathcal{M}$. Thus $\mathcal{M}$ is an order unit space.

By lemma 1 $\mathcal{A}$ is an order unit subspace of the order
unit space $\mathcal{M}$.  $\triangleright$

\medskip

{\bf Proposition 4.} {\it Let $A$ be a C$^*$-algebra on a Hilbert
space $H$, $\{p_i\}$ be an infinite orthogonal set of projections
in $A$ with the LUB $1$ in $B(H)$. Then $\mathcal{A}=\{\{p_\xi
ap_\eta\}: a\in A\}$ is a C$^*$-algebra, where the operation of
multiplication of $\mathcal{A}$ is defined as follows
$$
\cdot :<\{p_\xi ap_\eta\},\{p_\xi bp_\eta\}>\to \{p_\xi
abp_\eta\},  \{p_\xi ap_\eta\},\{p_\xi bp_\eta\}\in \mathcal{A}.
$$
}

{\bf Proof.} By lemma 4 in \cite{2} the map
$$
\mathcal{I}: a\in A\to \{p_\xi ap_\eta\}\in \mathcal{A}
$$
is a one-to-one map. In this case
$$
\mathcal{I}(a)\mathcal{I}(b)=\mathcal{I}(ab)
$$
by the definition of the operation of multiplication in
proposition 4, and $\mathcal{I}(a)=\{p_\xi ap_\eta\}$,
$\mathcal{I}(b)=\{p_\xi bp_\eta\}$, $\mathcal{I}(ab)=\{p_\xi ab
p_\eta\}$. Hence, the operation, introduced in proposition 4 is
associative multiplication and the map $\mathcal{I}$ is an
isomorphism of the algebras $A$ and $\mathcal{A}$.

By proposition 3 the isomorphism $\mathcal{I}$ is isometrical.
Therefore $\mathcal{A}$ is a C$^*$-algebra with this operation.
$\triangleright$

\medskip

{\it Example 1.} Let $H$ be a Hilbert space, $\{q_i\}$ be a
maximal orthogonal set of minimal projections in $B(H)$. Then
$\sup_i q_i=1$ and by lemma 4 in \cite{2} and proposition 4 the
algebra $\mathcal{B(H)}=\{\{q_iaq_j\}: a\in B(H)\}$ can be
identified with $B(H)$ as C$^*$-algebras in the sense of the map
$$
\mathcal{I}: a\in B(H)\to \{q_iaq_j\}\in \mathcal{B(H)}.
$$
In this case associative multiplication in $\mathcal{B(H)}$ is
defined as follows
$$
\cdot :<\{q_i aq_j\},\{q_i bq_j\}>\to \{q_i abq_j\}, \{q_i
aq_j\},\{q_i bq_j\}\in \mathcal{B(H)}.
$$
Let $a$, $b\in B(H)$, $q_iaq_j=\lambda_{ij}q_{ij}$,
$q_ibq_j=\mu_{ij}q_{ij}$, where $\lambda_{ij}$, $\mu_{ij}\in {\bf
C}$, $q_i=q_{ij}q_{ij}^*$, $q_j=q_{ij}^*q_{ij}$, for all indexes
$i$ and $j$. Then this operation of multiplication coincides with
the following bilinear operation
$$
\cdot :<\{q_i aq_j\},\{q_i bq_j\}>\to \{\sum_\xi \lambda_{i\xi}
\mu_{\xi j}q_{ij}\},  \{q_i aq_j\},\{q_i bq_j\}\in \mathcal{B(H)}.
$$

\medskip

{\it Remark 1.} Let  $A$ be a C$^*$-algebra on a Hilbert space
$H$, $\{p_i\}$ be an infinite orthogonal set of projections in $A$
with the LUB $1$ in $B(H)$. Then by proposition 4
$\mathcal{A}=\{\{p_\xi ap_\eta\}: a\in A\}$ is a C$^*$-algebra. In
this case the operation of involution on the algebra $\mathcal{A}$
coincides with the map
$$
\{p_\xi ap_\eta\}^*=\{p_\xi a^*p_\eta\},\,\,a\in A.
$$
Indeed, the identification $\mathcal{A}\equiv A$ gives us
$a=\{p_\xi ap_\eta\}$ and $a^*=\{p_\xi a^*p_\eta\}$ for all $a\in
A$. Then $\{p_\xi ap_\eta\}^*=a^*=\{p_\xi a^*p_\eta\}$ for any
$a\in A$. Let $\mathcal{A}_{sa}=\{\{p_\xi ap_\eta\}: a\in
A_{sa}\}$. Then $\mathcal{A}=\mathcal{A}_{sa}+i\mathcal{A}_{sa}$.
Indeed, $\{p_\xi ap_\eta\}^*=a^*=a=\{p_\xi ap_\eta\}$ for any
$a\in A_{sa}$.

Let $\mathcal{N}=\{\{p_\xi ap_\eta\}: a\in B(H)\}$. By lemma 4 in
\cite{2} and by proposition 4 $\mathcal{N}\equiv B(H)$. Therefore
we will assume that $\mathcal{N}=B(H)$. Let
$\mathcal{N}_{sa}=\{\{p_\xi ap_\eta\}: a\in B(H), \{p_\xi
ap_\eta\}^*=\{p_\xi ap_\eta\}\}$. Then
$\mathcal{N}=\mathcal{N}_{sa}+i\mathcal{N}_{sa}$. Note that
$\{p_\xi ap_\eta\}^*=\{p_\xi ap_\eta\}$ if and only if $(p_\xi
ap_\eta)^*=p_\eta ap_\xi$ for all $\xi$, $\eta$.

\medskip

{\bf Lemma 5.} {\it Let $H$ be a Hilbert space, $\{p_\xi\}$ be an
infinite orthogonal set of projections in $B(H)$ with the LUB 1.
Then associative multiplication of the algebra $\mathcal{N}$
(hence of the algebra $B(H)$) coincides with the operation
$$
\{p_\xi ap_\eta\}\star \{p_\xi bp_\eta\}=\{\sum_i p_\xi ap_i p_i
bp_\eta\},  \{p_\xi ap_\eta\},\{p_\xi bp_\eta\}\in \mathcal{N}
$$
where the sum $\sum$ in the right part of the equality is an
ultraweak limit of the net of finite sums of elements in the set
$\{p_\xi ap_i p_i bp_\eta\}_{\xi\eta}$.}

{\bf Proof.} Let $\{p_k\}_{k=1}^n$ be a finite subset of the set
$\{p_\xi\}$. Note that $\sup_i p_i=1$, i.e. the net of all finite
sums $\sum_{k=1}^n p_k$ of orthogonal projections in $\{p_\xi\}$
ultraweakly converges to the identity operator in $B(H)$. By the
ultraweakly continuity of the operator of multiplication $T(b)=ab,
b\in B(H)$, where $a\in B(H)$, the net of finite sums of elements
in $\{p_\xi ap_i p_i bp_\eta\}_{\xi\eta}$ ultraweakly converges
and $\sum_i p_\xi ap_i p_i bp_\eta=p_\xi abp_\eta$ for all $\xi$,
$\eta$. Hence the operation of multiplication $\star$ of the
algebra $\mathcal{N}$ coincides with the operation, introduced in
proposition 4. And the operation of associative multiplication,
introduced in proposition 4 coincides with multiplication in
$B(H)$ in the sense $\mathcal{N}\equiv B(H)$. $\triangleright$

\medskip

{\bf Proposition 6.} {\it Let $A$ be a C$^*$-algebra on a Hilbert
space $H$, $\{p_\xi\}$ be an infinite orthogonal set of
projections in $A$ with the LUB $1$ in $B(H)$. Then the operation
of associative multiplication of $\mathcal{A}$ coincides with the
operation of associative multiplication of $\mathcal{N}$ on $A$,
defined in lemma 5.}

{\bf Proof.} Let $\{p_\xi ap_\eta\}$, $\{p_\xi bp_\eta\}$ be
elements of $\mathcal{A}_{sa}$ and $\{p_k\}_{k=1}^n$ be a finite
subset of the set $\{p_\xi\}$ and $p=\sum_{k=1}^n p_k$. We have
the net of all finite sums $\sum_{k=1}^n p_k$ of orthogonal
projections in $\{p_\xi\}$ ultraweakly converges to the identity
operator in $B(H)$. Then for all $\xi$, $\eta$ the element
$\{p_\xi abp_\eta\}$ is an ultraweak limit in $B(H)$ of the net
$\{\sum_i p_\xi ap_i p_i bp_\eta\}$ of all finite sums
$\{\sum_{k=1}^n p_\xi ap_k p_k bp_\eta\}$ for all finite subsets
$\{p_k\}_{k=1}^n\subset \{p_\xi\}$, and the element $\{p_\xi
abp_\eta\}$ belongs to $\mathcal{A}$. Hence the assertion of
proposition 6 holds. $\triangleright$

\medskip

{\it Remark 2.} Let $A$ be a C$^*$-algebra on a Hilbert space $H$,
$\{p_i\}$ be an infinite orthogonal set of projections in $A$ with
the LUB $1$ in $B(H)$. Then by lemma 4 in \cite{2} the order and
the norm in the vector space $\sum_{i,j}^\oplus p_i Ap_j$ can be
introduced as follows: we write $\{a_{ij}\}\geq 0$, if this
element is zero or positive element in $B(H)$ in the sense of the
equality $B(H)=\sum_{\xi,\eta}^\oplus q_\xi B(H)q_\eta$, where
$\{q_\xi\}$ is an arbitrary maximal orthogonal set of minimal
projections in $B(H)$; $\Vert \{a_{ij}\}\Vert$ is equal to the
norm in $B(H)$ of this element in the sense of the equality
$B(H)=\sum_{\xi,\eta}^\oplus q_\xi B(H)q_\eta$ (example 1). By
lemmas 3 and 4 in \cite{2} they coincide with the order and the
norm defined in lemma 2 and proposition 3, respectively.

\medskip

{\it Remark 3.} Suppose that all conditions of remark 2 hold. Let
$\mathcal{B(H)}=\sum_{\xi,\eta}^\oplus q_\xi B(H)q_\eta$. Then
$B(H)\equiv\mathcal{B(H)}$, where $\mathcal{B(H)}=\{\{q_\xi
aq_\eta\}: a\in B(H)\}$. Also, we have $\sum_{ij}^\oplus p_i Ap_j$
is a Banach space and an order unit space (lemma 2, Proposition
3). Suppose that $\{q_\xi\}$ is a maximal orthogonal set of
minimal projections in $B(H)$ such that $p_i=\sup_\eta q_\eta$,
for some subset $\{q_\eta\}\subset \{q_\xi\}$, for all $i$. Note
that $B(H)\equiv\{\{p_iap_j\}: a\in B(H)\}=\sum_{ij}^\oplus p_i
B(H)p_j$. By propositions 4 and 6 the order unit space
$\mathcal{A}=\{\{p_iap_j\}: a\in A\}$ is closed with respect to
associative multiplication of $\sum_{ij}^\oplus p_i B(H)p_j$ (i.e.
$\mathcal{N}=\{\{p_iap_j\}: a\in B(H)\}$).

At the same time, the order unit space $\sum_{ij}^\oplus p_i Ap_j$
is the order unit subspace of the algebra $\sum_{ij}^\oplus p_i
B(H)p_j$.

Since $B(H)\equiv\sum_{ij}^\oplus p_i B(H)p_j$, we have
$\sum_{ij}^\oplus p_i B(H)p_j$ is a von Neumann algebra, and
without loss of generality, this algebra can be considered as the
algebra $B(H)$.

Note that if $\sum_{ij}^\oplus p_i Ap_j$ is closed with respect to
associative multiplication of $\sum_{ij}^\oplus p_i B(H)p_j$, then
$\sum_{ij}^\oplus p_i Ap_j$ is a C$^*$-algebra. Also, when we
consider the C$^*$-algebra $A$ with the conditions which are
listed above, then we have the algebra $\sum_{ij}^\oplus p_i
B(H)p_j$ (i.e. actually the algebra $B(H)$) and the vector space
$\sum_{ij}^\oplus p_i Ap_j$ as an order unit subspace of the
algebra $\sum_{ij}^\oplus p_i B(H)p_j$. Then we have
$$
\mathcal{A}\subseteq \sum_{ij}^\oplus p_i Ap_j\subseteq
\sum_{ij}^\oplus p_i B(H)p_j.
$$
Thus, further, when we say that $\sum_{ij}^\oplus p_i Ap_j$ is a
C$^*$-algebra we assume that the vector space $\sum_{ij}^\oplus
p_i Ap_j$ is closed with respect to associative multiplication of
the algebra $\sum_{ij}^\oplus p_i B(H)p_j$.

The involution in $\sum_{ij}^\oplus p_i B(H)p_j$ in the sense of
the identification $\sum_{ij}^\oplus p_i B(H)p_j\equiv B(H)$
coincides with the map
$$
\{a_{ij}\}^*=\{a_{ji}^*\}, \{a_{ij}\}\in \sum_{ij}^\oplus p_i
B(H)p_j.
$$
Indeed, there exists an element $a\in B(H)$ such that
$a=\{a_{ij}\}=\{p_iap_j\}$. Then $a^*=\{p_ia^*p_j\}$ in the sense
of $B(H)\equiv \mathcal{N}$. We have $a_{ij}=p_iap_j$,
$a_{ij}^*=p_ja^*p_i$ for all $i$, $j$. Therefore
$\{p_ia^*p_j\}=\{a_{ji}^*\}$. Hence $a^*=\{a_{ji}^*\}$. Let
$(\sum_{ij}^\oplus p_i B(H)p_j)_{sa}=\{\{a_{ij}\}: \{a_{ij}\}\in
\sum_{ij}^\oplus p_i B(H)p_j, \{a_{ij}\}^*=\{a_{ij}\}\}$. Then
$$
\sum_{ij}^\oplus p_i B(H)p_j=(\sum_{ij}^\oplus p_i
B(H)p_j)_{sa}+i(\sum_{ij}^\oplus p_i B(H)p_j)_{sa}.
$$

\medskip

{\bf Lemma 7.} {\it Let $A$ be a C$^*$-algebra on a Hilbert space
$H$, $\{p_i\}$ be an infinite orthogonal set of projections of $A$
with LUB $1$ in $B(H)$ and $(\sum_{ij}^\oplus p_i
Ap_j)_{sa}=\{\{a_{ij}\}: \{a_{ij}\}\in \sum_{ij}^\oplus p_i Ap_j,
\{a_{ij}\}^*=\{a_{ij}\}\}$. Then
$$
\sum_{ij}^\oplus p_i Ap_j=(\sum_{ij}^\oplus p_i
Ap_j)_{sa}+i(\sum_{ij}^\oplus p_i Ap_j)_{sa}. \,\,\,(1)
$$
In this case $\{a_{ij}\}^*=\{a_{ij}\}$ if and only if
$a_{ij}^*=a_{ji}$ for all $i$, $j$.}

{\bf Proof.} Let $\{a_{ij}\}\in \sum_{ij}^\oplus p_i Ap_j$. Since
$a_{ij}+a_{ji}\in A$, we have $a_{ij}+a_{ji}=a_1+ia_2$, where
$a_1$, $a_2\in (\sum_{ij}^\oplus p_i Ap_j)_{sa}$, for all $i$ and
$j$. Then
$a_{ij}+a_{ji}=p_ia_1p_j+p_ja_1p_i+i(p_ia_2p_j+p_ja_2p_i)$,
$a_1=p_ia_1p_j+p_ja_1p_i$, $a_2=p_ia_2p_j+p_ja_2p_i$ for all $i$
and $j$. Let $a^1_{ij}=p_ia_1p_j+p_ja_1p_i$,
$a^2_{ij}=p_ia_2p_j+p_ja_2p_i$ for all $i$ and $j$. Then by the
definition of $\sum_{ij}^\oplus p_i Ap_j$ we have $\{a^1_{ij}\}$,
$\{a^2_{ij}\}\in \sum_{ij}^\oplus p_i Ap_j$. In this case
$\{a^k_{ij}\}^*=\{a^k_{ij}\}$, $k=1,2$. Since the element
$\{a_{ij}\}\in \sum_{ij}^\oplus p_i Ap_j$ was chosen arbitrarily,
we have the equality (1).

The rest part of lemma 7 holds by the definition of the
self-adjoint elements $\{a^k_{ij}\}$, $k=1,2$. $\triangleright$

\medskip

{\bf Lemma 8.} {\it Let $H$ be a Hilbert space, $\{p_\xi\}$ be an
infinite orthogonal set of projections in $B(H)$ with the LUB $1$.
Then the operation of associative multiplication of the algebra
$\sum_{\xi,\eta}^\oplus p_\xi B(H)p_\eta$ (i.e. of the algebra
$B(H)$) coincides with the binary operation
$$
\cdot :<\{a_{\xi,\eta}\},\{b_{\xi,\eta}\}>\to \{\sum_i a_{\xi i}
b_{i\eta}\},  \{a_{\xi\eta}\},\{b_{\xi\eta} \}\in
(\sum_{\xi,\eta}^\oplus p_\xi B(H)p_\eta). \,\,\,\,\, (2)
$$
}

{\bf Proof.} Let $\{a_{\xi\eta}\},\{b_{\xi\eta} \}\in
(\sum_{\xi,\eta}^\oplus p_\xi B(H)p_\eta)$. We have
$$
B(H)\equiv \mathcal{N}\equiv \sum_{\xi,\eta}^\oplus p_\xi
B(H)p_\eta.
$$
Therefore, we can admit that
$B(H)=\mathcal{N}=\sum_{\xi,\eta}^\oplus p_\xi B(H)p_\eta$. There
exists elements $a$, $b$ in the algebra $B(H)$ such that $p_\xi
ap_\eta=a_{\xi\eta}$, $p_\xi bp_\eta=b_{\xi\eta}$ for all $\xi$,
$\eta$. Therefore $\{a_{\xi\eta}\}=\{p_\xi ap_\eta\}$,
$\{b_{\xi\eta}\}=\{p_\xi bp_\eta\}$. Then by lemma 5 associative
multiplication of the algebra $\sum_{\xi,\eta}^\oplus p_\xi
B(H)p_\eta$ (i.e. of the algebra $B(H)$) coincides with binary
operation (2). $\triangleright$

\medskip

{\bf Proposition 9.} (\cite{1}){\it Let $A$ be a von Neumann
algebra on a Hilbert space $H$, $\{p_i\}$ be an infinite
orthogonal set of projections of the algebra $A$ with LUB $1$.
Then $A=\sum_{\xi,\eta}^\oplus p_\xi Ap_\eta$.}

{\bf Proof.} Let $a$ be an element of $\sum_{\xi,\eta}^\oplus
p_\xi Ap_\eta$ and $a=\{a_{\xi\eta}\}$, where $a_{\xi\xi}=p_\xi
ap_\xi$, $a_{\xi\eta}=p_\xi ap_\eta$ for all $\xi$, $\eta$. We
have $a\in B(H)=\sum_{\xi,\eta}^\oplus p_\xi B(H)p_\eta$ and
$(\sum_{k=1}^n p_k)a(\sum_{k=1}^n p_k)\in A$ for any
$\{p_k\}_{k=1}^n\subset \{p_\xi\}$. Let
$$
b_n^\alpha=\sum_{kl=1}^n p_k^\alpha ap_l^\alpha=(\sum_{kl=1}^n
p_k^\alpha)a(\sum_{kl=1}^n p_k^\alpha)
$$
for all natural numbers $n$ and finite subsets $\{p_k^\alpha\}_{k=1}^n\subset
\{p_i\}$. Then by the proof of lemma 3 in \cite{2} the net
$(b_n^\alpha)$ ultraweakly converges to $a$ in $B(H)$. At the same
time $A$ is ultraweakly closed in $B(H)$. Therefore $a\in A$ and
$\sum_{\xi,\eta}^\oplus p_\xi Ap_\eta\subseteq A$.
$\triangleright$

{\bf Lemma 10.} {\it Let $A$ be a C$^*$-algebra on a Hilbert space
$H$, $\{p_\xi\}$ be an infinite orthogonal set of projections in
$A$ with the LUB $1$ in $B(H)$. Then, if projections of the set
$\{p_\xi\}$ are pairwise equivalent and for every index $\xi$ the
component $p_\xi Ap_\xi$ is a von Neumann algebra, then the vector
space $\sum_{\xi,\eta}^\oplus p_\xi Ap_\eta$ is closed with
respect to multiplication of the algebra $\sum_{\xi,\eta}^\oplus
p_\xi B(H)p_\eta$ and $\sum_{\xi,\eta}^\oplus p_\xi Ap_\eta$ is a
C$^*$-algebra.}

{\bf Proof.} First, note that $(p_\xi+p_\eta)A(p_\xi+p_\eta)$ is a
von Neumann algebra. Indeed, for any net $(a_\alpha)$ in $p_\xi
Ap_\eta$, weakly converging in $B(H)$ the net $(a_\alpha
x_{\xi\eta}^*)$ belongs to $p_\xi Ap_\xi$, where $x_{\xi\eta}$ is
an isometry in $A$ such that $x_{\xi\eta}x_{\xi\eta}^*=p_\xi$,
$x_{\xi\eta}^*x_{\xi\eta}=p_\eta$. Then, since the net $(a_\alpha
x_{\xi\eta}^*)$ weakly converges in $B(H)$, we have the weak limit
$b$ in $B(H)$ of the net $(a_\alpha x_{\xi\eta}^*)$ belongs to
$p_\xi Ap_\xi$. Hence $bx_{\xi\eta}\in p_\xi Ap_\eta$. It is easy
to see that $bx_{\xi\eta}$ is a weak limit in $B(H)$ of the net
$(a_\alpha)$. Hence $p_\xi Ap_\eta$ is weakly closed in $B(H)$.

Let $\{a_{\xi\eta}\}$, $\{b_{\xi\eta}\}\in (\sum_{\xi\eta}^\oplus
p_\xi Ap_\eta)$. We have
$$
\sum_{\xi\eta}^\oplus p_\xi Ap_\eta\subseteq \sum_{\xi\eta}^\oplus
p_\xi B(H)p_\eta=B(H).
$$
Therefore there exist elements $a$, $b$ in $\sum_{\xi\eta}^\oplus
p_\xi B(H)p_\eta$ (i.e. in the algebra $B(H)$) such that $p_\xi
ap_\eta=a_{\xi\eta}$, $p_\xi bp_\eta=b_{\xi\eta}$ for all $\xi$,
$\eta$. Therefore $\{a_{\xi\eta}\}=\{p_\xi ap_\eta\}$,
$\{b_{\xi\eta}\}=\{p_\xi bp_\eta\}$. We have
$$
\sum_i a_{\xi i} b_{i\eta}=p_\xi abp_\eta,
$$
calculated in $\sum_{\xi\eta}^\oplus p_\xi B(H)p_\eta$, belongs to
$p_\xi Ap_\eta$. Since the indices $\xi$, $\eta$ were chosen
arbitrarily and the product $\{p_\xi ap_\eta\}\{p_\xi
bp_\eta\}=ab$ belongs to $\sum_{\xi\eta}^\oplus p_\xi B(H)p_\eta$,
we have the product of the elements $a$ and $b$ belongs to
$\sum_{\xi,\eta}^\oplus p_\xi Ap_\eta$. Therefore
$\sum_{\xi\eta}^\oplus p_\xi Ap_\eta$ is closed with respect to
associative multiplication of the algebra $\sum_{\xi\eta}^\oplus
p_\xi B(H)p_\eta$. At the same time, $\sum_{\xi\eta}^\oplus p_\xi
Ap_\eta$ is a norm closed subspace of the algebra
$\sum_{\xi\eta}^\oplus p_\xi B(H)p_\eta=B(H)$. Hence
$\sum_{\xi\eta}^\oplus p_\xi Ap_\eta$ is a C$^*$-algebra and the
operation of multiplication in $\sum_{\xi\eta}^\oplus p_\xi
Ap_\eta$ can be defined as in lemma 8. $\triangleright$

\medskip

{\bf Theorem 11.} {\it Let $A$ be a C$^*$-algebra on a Hilbert
space $H$, $\{p_\xi\}$ be an infinite orthogonal set of
projections in $A$ with the LUB $1$ in $B(H)$. Then the following
statements hold:

1) Suppose that projections of the set $\{p_\xi\}$ are pairwise
equivalent and for any $\xi$ $p_\xi Ap_\xi$ is a von Nemann
algebra. Then $\sum_{\xi,\eta}^\oplus p_\xi Ap_\eta$ is a von
Neumann algebra,

2) if $\sum_{\xi,\eta}^\oplus p_\xi Ap_\eta$ is monotone complete
in $B(H)$ then $\sum_{\xi,\eta}^\oplus p_\xi Ap_\eta$ is a
C$^*$-algebra.}

{\bf Proof.} 1) Let $\{x_{\xi\eta}\}$ be a set of isometries in
$A$ such that $p_\xi=x_{\xi\eta}x_{\xi\eta}^*$,
$p_\eta=x_{\xi\eta}^*x_{\xi\eta}$ for all $\xi$, $\eta$. Let
$\xi$, $\eta$ be arbitrary indices. We prove that $p_\xi Ap_\eta$
is weakly closed. We have $p_\xi Ap_\eta p_\eta Ap_\xi \subseteq
p_\xi Ap_\xi$ and $p_\xi Ap_\eta =x_{\xi\eta}Ax_{\xi\eta}$. Let
$(a_\alpha)$ be a net in $p_\xi Ap_\eta$, weakly converging to an
element $a$ in $B(H)$. Then the exists a net $(b_\alpha)$ in
$p_\xi Ap_\eta$ such that $a_\alpha=x_{\xi\eta}b_\alpha
x_{\xi\eta}$ for all $\alpha$. By separately weakly continuity of
multiplication the net $(a_\alpha x_{\xi\eta}^*)$ weakly converges
to $ax_{\xi\eta}$ in $B(H)$. Since $(a_\alpha x_{\xi\eta}*)\subset
p_\xi Ap_\xi$ and $p_\xi Ap_\xi$ is weakly closed in $B(H)$, we
have $ax_{\xi\eta}^*\in p_\xi Ap_\xi$. Hence there exists an
element $b\in A$ such that
$ax_{\xi\eta}^*=x_{\xi\eta}bx_{\xi\eta}x_{\xi\eta}^*$. Then
$ax_{\xi\eta}^*x_{\xi\eta}=x_{\xi\eta}bx_{\xi\eta}x_{\xi\eta}^*x_{\xi\eta}
=x_{\xi\eta}bx_{\xi\eta}p_\eta=x_{\xi\eta}bx_{\xi\eta}\in p_\xi
Ap_\eta$. At the same time $a_\alpha p_\eta=a_\alpha$ for all
$\alpha$. Hence, $ap_\eta=a$. Since
$a=ax_{\xi\eta}^*x_{\xi\eta}=x_{\xi\eta}bx_{\xi\eta}\in p_\xi
Ap_\eta$, we have $a\in p_\xi Ap_\eta$. Also, since the net
$(a_\alpha)$ is chosen arbitrarily, we obtain the component $p_\xi
Ap_\eta$ is weakly closed in $B(H)$. Let $(a_\alpha)$ be a net in
$\sum_{\xi,\eta}^\oplus p_\xi Ap_\eta$, weakly converging to an
element $a$ in $B(H)$. Then for all $\xi$ and $\eta$ the net
$(p_\xi a_\alpha p_\eta)$ weakly converges to $p_\xi ap_\eta$ in
$B(H)$. In this case, by the previous part of the proof $p_\xi
ap_\eta\in p_\xi Ap_\eta$ for all $\xi$, $\eta$. Note that $a\in
\sum_{\xi,\eta}^\oplus p_\xi B(H)p_\eta$. Hence $a\in
\sum_{\xi,\eta}^\oplus p_\xi Ap_\eta$. Since the net $(a_\alpha)$
is chosen arbitrarily, we see that $\sum_{\xi,\eta}^\oplus p_\xi
Ap_\eta$ is weakly closed in $\sum_{\xi,\eta}^\oplus p_\xi
B(H)p_\eta\equiv B(H)$. Therefore by lemma 10
$\sum_{\xi,\eta}^\oplus p_\xi Ap_\eta$ is a von Neumann algebra.

Item 2) is a straightforward consequent of 1). $\triangleright$

\medskip

{\bf Proposition 12.} {\it Let $A$ be a monotone complete
C$^*$-algebra on a Hilbert space $H$, $\{p_\xi\}$ be an infinite
orthogonal set of projections in $A$ with the LUB $1$ in $B(H)$.
Then the order unit space $\sum_{\xi,\eta}^\oplus p_\xi Ap_\eta$
is monotone complete.}

{\bf Proof.} It is clear that the C$^*$-subalgebra $p_\xi Ap_\xi$
is also monotone complete for any index $\xi$. Let
$\{p_k\}_{k=1}^n$ be a finite subset of $\{p_\xi\}$ and
$p=\sum_{k=1}^n p_k$. Then the C$^*$-subalgebra $pAp$ is also
monotone complete.

Let $(a_\alpha)$ be a bounded monotone increasing net in
$\sum_{\xi,\eta}^\oplus p_\xi Ap_\eta$. Since for any finite
subset $\{p_k\}_{k=1}^n\subseteq \{p_\xi\}$ the subalgebra
$(\sum_{k=1}^n p_k) A(\sum_{k=1}^n p_k)$ is monotone complete, we
have
$$
\sup_\alpha (\sum_{k=1}^n p_k)a_\alpha (\sum_{k=1}^n p_k)\in
(\sum_{k=1}^n p_k) A(\sum_{k=1}^n p_k).
$$
Hence, $\{a_{\xi\eta}\}=\{\sup_\alpha p_\xi a_\alpha p_\xi\}\cup
\{p_\xi(\sup_\alpha (p_\xi+p_\eta) a_\alpha
(p_\xi+p_\eta))p_\eta\}_{\xi\neq\eta}$ is an element of the order
unit space $\sum_{\xi,\eta}^\oplus p_\xi Ap_\eta$. It can be
checked straightforwardly using the definition of the order in
$\sum_{\xi,\eta}^\oplus p_\xi Ap_\eta$ that the element
$\{a_{\xi\eta}\}$ is LUB of the net $(a_\alpha)$. Since the net
$(a_\alpha)$ was chosen arbitrarily, we obtain the order unit
space $\sum_{\xi,\eta}^\oplus p_\xi Ap_\eta$ is monotone compete.
$\triangleright$

\medskip

{\bf Theorem 13.} {\it Let $A$ be a monotone complete
C$^*$-algebra on a Hilbert space $H$, $\{p_\xi\}$ be an infinite
orthogonal set of projections in $A$ with the LUB $1$ in $B(H)$.
Suppose that projections of the set $\{p_\xi\}$ are pairwise
equivalent and $A$ is not a von Nemann algebra. Then $A\neq
\sum_{\xi,\eta}^\oplus p_\xi Ap_\eta$ (i.e.
$\mathcal{A}:=\{\{p_\xi ap_\eta\}: a\in A\}\neq
\sum_{\xi,\eta}^\oplus p_\xi Ap_\eta$). }

{\bf Proof.} By the condition there exists a bounded monotone
increasing net $(a_\alpha)$ of elements in $A$, the LUB $\sup_{A}
a_\alpha$ in $A$ and the LUB $\sup_{\sum_{\xi\eta}^\oplus p_\xi
B(H)p_\eta} a_\alpha$ in $\sum_{\xi\eta}^\oplus p_\xi B(H)p_\eta$
of which are different. Otherwise $A$ is a von Nemann algebra.

By the definition of the order in $\sum_{\xi\eta}^\oplus p_\xi
B(H)p_\eta$ there exists a projection $p\in \{p_\xi\}$ such that
the LUB $\sup_{pAp} pa_\alpha p$ in $pAp$ and the LUB
$\sup_{pB(H)p} p a_\alpha p$ in $pB(H)p$ of the bounded monotone
increasing net $(pa_\alpha p)$ of elements in $pAp$ are different.
Indeed, let $a=\sup_A a_\alpha$, $b=\sup_{\sum_{\xi\eta}^\oplus
p_\xi B(H)p_\eta} a_\alpha$. Since $A\subseteq
\sum_{\xi\eta}^\oplus p_\xi B(H)p_\eta$, we have $b\leq a$ and
$0\leq a-b$. Hence, if $p_\xi (a-b) p_\xi=0$ for all $\xi$, then
$p_\xi (a-b)=(a-b)p_\xi=0$. Therefore by lemma 2 in \cite{2}
$a-b=0$, i.e. $a=b$. Hence $pAp$ is not a von Nemann algebra.

There exists an infinite orthogonal set $\{e_i\}$ of projections
in $pAp$, the LUB $\sup_{pAp} e_i$ in $pAp$ and the LUB
$\sup_{pB(H)p} e_i$ in $pB(H)p$ of which are different. Otherwise
$pAp$ is a von Neumann algebra.

Indeed, any maximal commutative subalgebra $A_o$ of $pAp$ is
monotone complete. For any normal positive linear functional
$\rho\in B(H)$ and for any infinite orthogonal set $\{q_i\}$ of
projections in $A_o$ we have $\rho(\sup_i q_i)=\sum_i \rho(q_i)$,
where $\sup_i q_i$ is the LUB of the set $\{q_i\}$ in $A_o$. Hence
by the theorem on extension of a $\sigma$-additive measure to a
normal linear functional $\rho\vert_{A_o}$ is a normal functional
on $A_o$. Hence $A_o$ is a commutative von Neumann algebra. At the
same time the maximal commutative subalgebra $A_o$ of the algebra
$pAp$ is chosen arbitrarily. Therefore by \cite{3} $pAp$ is a von
Neumann algebra. What is impossible.

Let $\{x_{\xi\eta}\}$ be a set of isometries in $A$ such that
$p_\xi=x_{\xi\eta}x_{\xi\eta}^*$,
$p_\eta=x_{\xi\eta}^*x_{\xi\eta}$ for all $\xi$, $\eta$, and
$p_1=p$. Let $\{x_{1\xi}\}$ be the subset of the set
$\{x_{\xi\eta}\}$ such that $p_1=x_{1\xi}x_{1\xi}^*$,
$p_\xi=x_{1\xi}^*x_{1\xi}$ for all $\xi$. Without loss of
generality we regard that the set of indices $i$ for $\{e_i\}$ is
a subset of the set of indices $\xi$ for $\{p_\xi\}$. Let
$\{e_ix_{1i}\}$ be a set of all components of some infinite
dimensional matrix $\{a_{\xi\eta}\}$, the components, which are
not presented, are zeros and $\{x_{1i}^*e_i^*\}$ be also a similar
matrix, which  coincides with $\{a_{\xi\eta}^*\}$. We have $\sum_i
e_ix_{1i}x_{1i}^*e_i^*=\sum_i e_ip_1e_i^*=\sum_i e_ie_i^*=\sum_i
e_i\leq \sup_{pAp} e_i$. Therefore $\{a_{\xi\eta}\}\in
\sum_{\xi,\eta}^\oplus p_\xi Ap_\eta$. Then $\{a_{\xi\eta}^*\}\in
\sum_{\xi,\eta}^\oplus p_\xi Ap_\eta$. Therefore if
$\{a_{\xi\eta}\}\in A$ (i.e. in $\mathcal{A}:=\{\{p_\xi ap_\eta\}:
a\in A\}$), then the product
$\{a_{\xi\eta}\}\cdot\{a_{\xi\eta}^*\}$ in $\sum_{ij}^\oplus p_i
B(H)p_j$ belongs to $\sum_{\xi,\eta}^\oplus p_\xi Ap_\eta$. In
this case the infinite dimensional matrix
$\{a_{\xi\eta}\}\cdot\{a_{\xi\eta}^*\}$ contains the component
$\sum_i e_ix_{1i}\cdot x_{1i}^*e_i^*$ such that $\sum_i
e_ix_{1i}\cdot x_{1i}^*e_i^*=p_1(\sum_i e_ix_{1i}\cdot
x_{1i}^*e_i^*)p_1$. Consequently,
$p_1(\{a_{\xi\eta}\}\cdot\{a_{\xi\eta}^*\})p_1=\sum_i
e_ix_{1i}\cdot x_{1i}^*e_i^*$. Hence $\sum_i e_ix_{1i}\cdot
x_{1i}^*e_i^*\in p_1 (\sum_{\xi,\eta}^\oplus p_\xi
Ap_\eta)p_1=p_1Ap_1$. Since $\sum_i e_ix_{1i}\cdot
x_{1i}^*e_i^*=\sum_i e_ip_1e_i^*=\sum_i e_ie_i^*=\sum_i e_i$, we
obtain $\sum_i e_i\in p_1Ap_1$, i.e. $\sup_{pB(H)p} e_i\in
p_1Ap_1$. The last statement is a contradiction. Therefore
$\{a_{\xi\eta}\}\notin A$. Hence $A\neq \sum_{\xi,\eta}^\oplus
p_\xi Ap_\eta$ (i.e. $\mathcal{A}:=\{\{p_\xi ap_\eta\}: a\in
A\}\neq \sum_{\xi,\eta}^\oplus p_\xi Ap_\eta$).  $\triangleright$

\medskip

The following assertion follows by theorem 13 and it's proof.

{\bf Corollary 14.} {\it Let $A$ be a C$^*$-algebra on a Hilbert
space $H$, $\{p_\xi\}$ be an infinite orthogonal set of
projections in $A$ with the LUB $1$ in $B(H)$. Then the following
statements hold:

1) suppose that the order unit space $\sum_{\xi,\eta}^\oplus p_\xi
Ap_\eta$ is monotone complete and there exists a bounded monotone
increasing net $(a_\alpha)$ of elements in $\sum_{\xi,\eta}^\oplus
p_\xi Ap_\eta$, the LUB $\sup_{\sum_{\xi,\eta}^\oplus p_\xi
Ap_\eta} a_\alpha$ in $\sum_{\xi,\eta}^\oplus p_\xi Ap_\eta$ and
the LUB $\sup_{\sum_{\xi\eta}^\oplus p_\xi B(H)p_\eta} a_\alpha$
in $\sum_{\xi\eta}^\oplus p_\xi B(H)p_\eta$ of which are
different. Then the vector space $\sum_{\xi,\eta}^\oplus p_\xi
Ap_\eta$ is not closed with respect to multiplication of the
algebra $\sum_{\xi,\eta}^\oplus p_\xi B(H)p_\eta$,

2) if $\sum_{\xi,\eta}^\oplus p_\xi Ap_\eta$ is a C$^*$-algebra
then this algebra is a von Neumann algebra.}

\bigskip

\section{Application}

\medskip

Let $n$ be an infinite cardinal number, $\Xi$ a set of indices of
cardinality $n$. Let $\{e_{ij}\}$ be a set of matrix units such
that $e_{ij}$ is a $n\times n$-dimensional matrix, i.e.
$e_{ij}=(a_{\alpha\beta})_{\alpha\beta\i\Xi}$, $(i,j)$-s component
of which is $1$, i.e. $a_{ij}=1$, and the rest components are
zero. Let $X$ be a hyperstonean compact, $C(X)$ the commutative
algebra of all complex-valued continuous functions on the compact
$X$ and
$$
\mathcal{M}=\{\sum_{ij\in\Xi}\lambda_{ij}(x)e_{ij}: (\forall ij
\lambda_{ij}(x)\in C(X))
$$
$$
(\exists K\in R)(\forall m\in N)(\forall
\{e_{kl}\}_{kl=1}^m\subseteq \{e_{ij}\})\Vert\sum_{kl=1\dots
m}\lambda_{kl}(x)e_{kl}\Vert\leq K\},
$$
where $\Vert\sum_{kl=1\dots m}\lambda_{kl}(x)e_{kl}\Vert\leq K$
means $(\forall x_o\in X) \Vert \sum_{kl=1\dots
m}\lambda_{kl}(x_o)e_{kl}\Vert\leq K$. The set $\mathcal{M}$ is a
vector space with pointwise algebraic operations. The map
$\Vert\,\,\, \Vert : \mathcal{M}\to {\bf R}_+$ defined as
$$
\Vert a \Vert = \sup_{\{e_{kl}\}_{kl=1}^n\subseteq
\{e_{ij}\}}\Vert\sum_{kl=1}^n \lambda_{kl}(x)e_{kl}\Vert,
$$
is a norm on the vector space $\mathcal{M}$, where $a\in
\mathcal{M}$ and $a=\sum_{ij\in\Xi}\lambda_{ij}(x)e_{ij}$.

{\bf Theorem 15.} {\it $\mathcal{M}$ is a von Neumann algebra of
type I$_n$ and $\mathcal{M}=C(X)\otimes M_n({\bf C})$.}

{\it Proof.} It is easy to see that the set  $\mathcal{M}$ is a
vector space with the componentwise algebraic operations. It is
known that the vector space $C(X,M_n({\bf C}))$ of continuous
matrix-valued maps on the compact $X$ is a C$^*$-algebra. Let
$A=C(X,M_n({\bf C}))$ and $e_i$  be a $e_{ii}$-valued constant map
on $X$, i.e. $e_i$ is a projection of the algebra $A$. A
C$^*$-algebra $A$ can be embedded in $B(H)$ for some Hilbert space
$H$ such that $\{e_i\}$ is an orthogonal set of projections with
$\sup_i e_i=1$ in $B(H)$. Then $\sum_{ij}^\oplus
e_iAe_j=\mathcal{M}$ and $\sum_{ij}^\oplus e_iAe_j$ can be
embedded in $B(H)$. For any $i$ $e_iAe_i=C(X)e_i$, i.e. the
component $e_iAe_i$ is weakly closed in $B(H)$. Hence, by theorem
11 the image of vector space $\mathcal{M}$ in $B(H)$ is a von
Neumann algebra. Hence, $\mathcal{M}$ is a von Neumann algebra.
Note, that the set $\{e_i\}$ is a maximal orthogonal set of
abelian projections with central support 1. Hence, $\mathcal{M}$
is a von Neumann algebra of type I$_n$. Moreover the center
$Z(\mathcal{M})$ of the algebra $\mathcal{M}$ is isomorphic to
$C(X)$ and $\mathcal{M}=C(X)\otimes M_n({\bf C})$.
$\triangleright$

\bigskip

Andizhan State University, Andizhan, Uzbekistan

arzikulovfn@rambler.ru, arzikulovfn@yahoo.com

\end{document}